\documentclass{amsart}

\usepackage{amssymb}
\usepackage[all]{xy}

\newtheorem{thm}{Theorem}

\newtheorem{lem}[thm]{Lemma}
\newtheorem{cor}[thm]{Corollary}

\newtheorem{prop}[thm]{Proposition}

\theoremstyle{definition}
\newtheorem{defn}[thm]{Definition}

\newtheorem{say}[thm]{}
\newtheorem{exmp}[thm]{Example}

\newtheorem{ques}[thm]{Question}    

\newtheorem{rem}[thm]{Remark}          

\newtheorem{notation}[thm]{Notation}   
  
\newtheorem{defn-thm}[thm]{Definition--Theorem}  
\newtheorem{defn-lem}[thm]{Definition--Lemma}  

\theoremstyle{remark}

\renewcommand{\o}[0]{{\mathcal O}} 

\renewcommand{\r}[0]{{\mathbb R}}

\newcommand{\p}[0]{{\mathbb P}}

\newcommand{\q}[0]{{\mathbb Q}}
\newcommand{\map}[0]{\dasharrow}
\newcommand{\qtq}[1]{\quad\mbox{#1}\quad}

\newcommand{\supp}[0]{\operatorname{Supp}}    
    
\newcommand{\codim}[0]{\operatorname{codim}}

\newcommand{\cent}[0]{\operatorname{center}}

\newcommand{\ex}[0]{\operatorname{Ex}}    
\newcommand{\diff}[0]{\operatorname{Diff}}

\newcommand{\simq}[0]{\sim_{\q}}

\newcommand{\coeff}[0]{\operatorname{coeff}}

\def\into{\DOTSB\lhook\joinrel\to}

\def\loccoh#1.#2.#3.#4.{H^{#1}_{#2}(#3,#4)}

\DeclareMathAlphabet{\mathchanc}{OT1}{pzc}%
                                {m}{it}

\usepackage[all]{xy}\xyoption{dvips}

\newcommand{\dex}[0]{\operatorname{DEx}}

\begin{document}
\bibliographystyle{amsalpha}

 \title{Minimal models of semi-log-canonical pairs}
 \author{Florin Ambro and J\'anos Koll\'ar}

\begin{abstract} We compare the minimal model of a log canonical pair with the minimal model of its reduced boundry. These results are then used 
to study the existence of the minimal model of a semi-log-canonical pair using its normalization.  
\end{abstract}

 \maketitle

In birational geometry, it is frequently necessary to work not just with
log canonical pairs  $(X, \Delta)$, but with their non-normal variants, called
{\it semi-log-canonical pairs.} Such pairs appear when one tries to compactify the moduli spaces of varieties and in  inductive arguments.

Many properties of log canonical pairs have been generalized to 
the semi-log-canonical setting  \cite{ambro, amb-surv, MR3238112, fujinobook, kk-singbook}, but it was observed in 
\cite{MR2824967} that log canonical rings of semi-log-canonical pairs are not always  finitely generated and some flips of semi-log-canonical pairs do not exist. Note that, by contrast,  abundance holds for a  semi-log-canonical pair iff it holds for its normalization; this was proved in increasing generality in
\cite{MR966580, k-etal, MR1284817, fujino-sl-abund, MR3048544, fuj-gon-cbf, hacon-xu-brep}.

The aim of this note is to describe some conditions that guarantee the existence of minimal models  for certain semi-log-canonical pairs.  
Our assumptions are  rather restrictive, but they may be close to being optimal.
The key  is to understand an even simpler question involving log canonical pairs: {\it How does the boundary of a log canonical pair change under a flip?}

This is a very natural problem, that first appeared explicitly in Tsunoda's treatment of semi-stable flips \cite{MR946246}, later in Shokurov's approach
that reduces flips to special flips  \cite{shok-MR1162635, k-etal} and  in 
\cite[Sec.4]{hmx-bounded}; see also \cite{ber-pau}. 

We are thus led to the following general questions.

\begin{ques}\label{aq1}
 Let $(X, D+\Delta)$ be an  lc pair that is projective over a base scheme $S$  with relatively ample divisor  $H$, where  
 all divisors in $D$ appear with coeffcient 1. 
Set 
$(X^0, D^0+\Delta^0):=(X, D+\Delta)$ and for $i=1,\dots, m$ let 
$$
\phi^i: (X^{i-1}, D^{i-1}+\Delta^{i-1})\map (X^{i}, D^{i}+\Delta^{i})
$$
be the steps of the $(X, D+\Delta)$-MMP with scaling of $H$; see Definition \ref{say1}. 
Let  $\rho:\bar D\to D$ be the normalization.
Do the restrictions
$$
\phi^i_{ D}:= \phi^i|_{\bar D^{i-1}}:(\bar D^{i-1},\diff_{\bar D}\Delta^{i-1})\map 
(\bar D^{i},\diff_{\bar D}\Delta^{i})
$$
form  the steps of the MMP  
starting with $(\bar D^{0},\diff_{\bar D}\Delta^{0}):=(\bar D,\diff_{\bar D}\Delta)$ and with scaling of $\rho^*H$?
\end{ques}

\begin{notation} \label{rem1}
We follow the terminology and notation of \cite{km-book, kk-singbook}.

From now on, whenever we write a divisor as $D+\Delta$, we assume that all
irreducible components of $D$  appear with coefficient 1 ($\Delta$ may also contain divisors   with coefficient 1).

Let   $\rho:\bar D\to D$ denote the normalization. 
The {\it different} of $\Delta$ on $\bar D$ is denoted by $\diff_{\bar D}\Delta$. It  is a $\q$-divisor on $\bar D$  that satisfies a natural
$\q$-linear equivalence
$$
K_{\bar D}+\diff_{\bar D}\Delta\simq \rho^*\bigl(K_X+D+\Delta\bigr).
\eqno{(\ref{rem1}.1)}
$$
See \cite[4.2]{kk-singbook} for a precise definition and its main properties.
In order to avoid secondary sub and superscripts, we usually write
$\diff_{\bar D}\Delta^{i}$ instead of the more precise
$\diff_{\bar D^i}\Delta^{i}$.

In the original definition, a step of the MMP corresponds to an extremal ray \cite{ckm}. By (\ref{rem1}.1), 
any contraction of an extremal ray on $X$ induces the contraction of an extremal face on $D$, but the face may well have dimension $>1$.  In an MMP with scaling of an ample divisor, the steps  correspond to certain contractions of  extremal faces. The divisor $H$  plays a very minor role in the sequel, but it makes it possible for us to tell exactly which MMP steps we get. 
\end{notation}

It turns out that a positive answer to  Question \ref{aq1} gives a positive answer to
the following problem on slc pairs.

\begin{ques} \label{aq2}
Let $(X, \Delta)$ be an slc pair  that is projective over a base scheme $S$
with normalization  $\pi: (\bar X, \bar D+\bar \Delta)\to (X, \Delta)$, conductor $\bar D\subset \bar X$  and $H$ an ample divisor on $X$. Set
$(\bar X^0, \bar D^0+\bar \Delta^0):=(\bar X, \bar D+\bar \Delta)$ and
for $i=1,\dots, m$ let 
$$
\bar\phi^i: (\bar X^{i-1}, \bar D^{i-1}+\bar \Delta^{i-1})\map (\bar X^{i}, \bar D^{i}+\bar \Delta^{i})
$$
be the steps of the $(\bar X, \bar D+\bar \Delta)$-MMP with scaling of $\pi^*H$.
Do we get 
$$
\phi^i: (X^{i-1}, \Delta^{i-1})\map (X^{i}, \Delta^{i}),
$$
which form the  steps of the $(X,\Delta)$-MMP  
 with scaling of $H$ and such that
$(\bar X^{i}, \bar D^{i}+\bar \Delta^{i}) $ is the normalization of $(X^{i}, \Delta^{i}) $?
\end{ques}

\begin{exmp}\label{xmp1}
We give 2 types of examples showing that in  Question \ref{aq1}
we usually do not get the steps of the $(\bar D,\diff_{\bar D}\Delta)$-MMP.
\medskip

(\ref{xmp1}.1) Start with a smooth variety $X'$, a smooth divisor $D'\subset X'$
and another  smooth divisor $C'\subset D'$. Assume that $K_{X'}+D'$ 
 is ample.
Set $X:=B_{C'}X'$ with exceptional divisor $E$ and let
$D\subset X$ denote the birational transform of $D'$.

For any $1\geq \epsilon>0$,
 $(X, D+\epsilon E)$ is an lc pair whose canonical model is $(X', D')$
and $(D', 0)$ is its own
canonical model.

However,  $(D, \epsilon\diff_DE)\cong (D', \epsilon C')$ is different from $(D', 0)$.

Note further that $K_{X}+D$ is the pull-back of $K_{X'}+D'$, hence
semiample and big. Thus the stable base locus of
$K_{X}+D+\epsilon E$ is  $E$. If $1> \epsilon>0$ then the only
log canonical center of $(X, D+\epsilon E)$ is $D$ and the other
log centers are  $E$ and $E\cap D$; see Definition \ref{l.l.cent.defn}. 
Thus the stable base locus contains the log centers but not the 
log canonical center.

Here are some  concrete examples.  

(\ref{xmp1}.1.1) Let $X'$ be a smooth surface,  $D'\subset X'$ a smooth rational curve and $C'\subset D'$ a set of 3 points.
Then  $(D, \diff_DE)\cong (D',  C')$ has ample log canonical class 
but $(D', 0)\cong (\p^1, 0)$ has negative log canonical class.

(\ref{xmp1}.1.2) For $\dim X'\geq 3$  it can also happen that the $(D, \epsilon\diff_DE)$-MMP tells us to contract $C'$.
Take $X'=\p^3$ and let $D'\subset X'$ be a smooth surface of degree 5 that contains a line $C'$. Then the self-intersection of $C'$ is $-3$, thus for $1\geq \epsilon>\frac13$
the first (an only) step of the $(D, \epsilon\diff_DE)$-MMP is to contract $C'$.
\medskip

(\ref{xmp1}.2) Let $B$ be  a smooth curve and $f:X\to B$ be a flat family of surfaces with quotient singularities and such that $K_X$ is $\q$-Cartier.

Let  $g:X\to Z$ be a flipping contraction.
(For  concrete examples, see \cite[2.7]{km-book} or the
 list in \cite{km-flips}.)
Thus there is a closed point $0\in B$ such that
$g$ is an isomorphism over $B\setminus \{0\}$. 
Set $D:=X_0$ and let $C\subset D$ denote the flipping curve.
Our example is the pair $(X, D)$. Here $\diff_D0=0$, hence we need  to compare the MMP for $(X, D)$ with the MMP for $(D, 0)$. 

Over $0\in B$ we have a birational
contraction $g_0:X_0\to Z_0$ that contracts 
$C\subset X_0$ to a point. Moreover
$(C\cdot K_{X_0})=(C\cdot K_X)<0$, thus
$Z_0$ is again log terminal and the contraction
$g_0:X_0\to Z_0$ is a step in the MMP for $X_0=D$.

However, since  $g:X\to Z$ a flipping contraction,
the special fiber of the flip  $g^+:X^+\to Z$
is another surface  $X_0^+\to Z_0$ with a new exceptional curve
$C^+\subset X_0^+$ such that 
$\bigl(C^+\cdot K_{X^+_0})=(C^+\cdot K_{X^+})>0$. Thus $X_0^+$ is not the 
canonical model of $X_0$ and $X_0\map X_0^+$ is not even a step of
any minimal model program.

We can easily arrange that $K_{X^+}$ is ample. In this case the stable base locus of $K_X$ is the flipping curve $C\subset X_0=D$. 
 The only
log canonical center of $(X, D)$ is $D$ which is
not contained in the  stable base locus of $K_X$.

It is easy to see that $D$ must have at least 1 non-canonical singularity that is also contained in $C$.  This gives a 0-dimensional log center of $(X, D)$
that is contained in the  stable base locus.
\end{exmp}

\begin{exmp} Every counter example  to Question \ref{aq1},
where $D$ is normal, gives a  counter example  to Question \ref{aq2} as follows.

Let $b\in B$ be a smooth, projective, pointed  curve of genus $\geq 1$. 
We can glue  $(X, D+\Delta)$ to  $(B\times D, \{b\}\times D+B\times \diff_D\Delta)$ along $D$ to get an slc pair 
$(Y, \Delta_Y)$ 
whose nomalization is
the disjoint union of $(X, D+\Delta)$ and  $(B\times D, \{b\}\times D+B\times \diff_D\Delta)$. On  $(X, D+\Delta)$ we get the steps of the  $(X, D+\Delta)$-MMP 
$$
\phi^i: (X^{i-1}, D^{i-1}+\Delta^{i-1})\map (X^{i}, D^{i}+\Delta^{i})
$$
and these restrict to
$$
\phi^i_{ D}: ( D^{i-1},\diff_{ D}\Delta^{i-1})\map 
( D^{i},\diff_{ D}\Delta^{i}).
$$
Let us denote the  steps of the $(D, \diff_D\Delta)$-MMP by
$$
\psi_i: ( D_{i-1},\diff_{ D}\Delta_{i-1})\map 
( D_{i},\diff_{ D}\Delta_{i}).
$$
Then the steps of the $(B\times D, \{b\}\times D+B\times \diff_D\Delta)$-MMP are given by
$$
\bigl(B\times D_{i-1}, \{b\}\times D_{i-1}+B\times \diff_D\Delta_{i-1}\bigr)\map
\bigl(B\times D_i, \{b\}\times D_i+B\times \diff_D\Delta_i\bigr).
$$
If 
$( D^{i},\diff_{ D}\Delta^{i})\ncong ( D_{i},\diff_{ D}\Delta_{i})$,
then we  can not glue the resulting
pairs  
$$
 (X^{i}, D^{i}+ \Delta^{i})\qtq{and}
\bigl(B\times D_i, \{b\}\times D_i+B\times \diff_D\Delta_i\bigr).
$$
Thus the $(Y, \Delta_Y)$-MMP does not exist.  
\end{exmp}

We give positive answers to  Questions \ref{aq1} and \ref{aq2}
when the singularities of $(X, D+\Delta)$ (resp.\ of $(\bar X, \bar D+\bar\Delta)$) are mild along the exceptional locus of $\phi$ (resp.\ of $\bar \phi$).
We use discrepancies to make this assertion precise.

\begin{defn}\label{l.l.cent.defn}
 Let  $(X, \Theta)$  be an lc pair. An irreducible subvariety $W\subset X$ is called a {\it log canonical center}  (resp.\  a {\it log center}) of  $(X, \Theta)$ if there is a divisor $E$ over $X$ such that
$\cent_XE=W$ and 
$a(E, X, \Theta)=-1$   (resp.\ $a(E, X, \Theta)<0$).

Assume next that $\Theta=D+\Delta$ and let $\rho: \bar D\to D$ denote the normalization.
By adjunction \cite[4.9]{kk-singbook}, $W\subset \bar D$ is a log center of
$(\bar D,\diff_{\bar D}\Delta)$ iff $\rho(W)$ is a log center of  $(X, D+\Delta)$. 
See \cite[Chap.7]{kk-singbook} for more on log centers. 
\end{defn}

From now on we assume  that the base scheme $S$ is essentially of finite type over a field of characteristic 0. Our main result is the following.

\begin{thm}\label{at1}
Using the notation and assumptions of  Question \ref{aq1}, 
assume in addition that the intersection of $D$ with the exceptional locus of
$$
\Phi^m:=\phi^m\circ\cdots\circ \phi^1: X\map X^{m}
$$
does not contain any log center of  $(X, D+\Delta)$.
Then the maps 
$$
\phi^i_{\bar D}: (\bar D^{i-1},\diff_{\bar D}\Delta^{i-1})\map 
(\bar D^{i},\diff_{\bar D}\Delta^{i})
$$
form  the steps of the MMP  
starting with $(\bar D^{0},\diff_{\bar D}\Delta^{0}):=(\bar D,\diff_{\bar D}\Delta)$ and with scaling of $\rho^*H$.
\end{thm}

\begin{rem}
As the Examples (\ref{xmp1}.1.1--2) show, we need to avoid all log centers, not just the log canonical centers.

It can happen that  $\phi^i$ is an isomorphism along  $D^{i-1}$. 
Thus the precise claim is that each $\phi^i_{\bar D} $ is either an isomorphism or an MMP step.  (The literature is  somewhat inconsistent. Usual  definitions of MMP steps  allow isomorphisms, but 
in many statements they are tacitly excluded.)

\end{rem}

\begin{thm}\label{at2}
Using the notation and assumptions of  Question \ref{aq2}, 
assume in addition that the intersection of $\bar D$ with the exceptional locus of
$$
\Phi^m_{\bar X}:=\bar\phi^m\circ\cdots\circ \bar\phi^1:
 \bar X\map \bar X^{m}
$$
does not contain any log center of  $(\bar X,\bar D+\bar\Delta)$.

Then the first $m$ steps  of the $(X,\Delta)$-MMP with scaling of $H$ exist
$$
\phi^i: (X^{i-1}, \Delta^{i-1})\map (X^{i}, \Delta^{i}),
$$
and 
$ (\bar X^{i}, \bar D^{i}+\bar \Delta^{i})$ is the normalization of
$(X^{i}, \Delta^{i}) $. 
\end{thm}

Proof. Let $(X,  \Delta)$ be an slc pair with 
normalization $(\bar X, \bar D+\bar \Delta)\to (X,  \Delta)$,
where $\bar D\subset \bar X$ is the conductor. 
Let $\rho:\bar D^n\to \bar D$ denote its normalization.

The gluing theory of  \cite[Chap.5]{kk-singbook} says that 
there is a (regular) involution 
$$
\tau:(\bar D^n,\diff_{\bar D^n}\bar\Delta)\to 
(\bar D^n,\diff_{\bar D^n}\bar\Delta),
$$ 
and $X$ is obtained from $\bar X$ by identifying the  equivalence classes of the relation generated by $\tau$ on $\bar X$.

Next let 
$$
\bar\phi^i: (\bar X^{i-1}, \bar D^{i-1}+\bar \Delta^{i-1})\map (\bar X^{i}, \bar D^{i}+\bar \Delta^{i})
$$
be the steps of the $(\bar X, \bar D+\bar \Delta)$-MMP with scaling of $\pi^*H$ and assume that Theorem \ref{at1} applies. Then 
$$
\bar\phi^i_D: \bigl((\bar D^{i-1})^n, \diff_{\bar D^n}\bar \Delta^{i-1}\bigr)
\map 
 \bigl((\bar D^{i})^n, \diff_{\bar D^n}\bar \Delta^{i}\bigr)
$$
are steps of the $ \bigl(\bar D^n, \diff_{\bar D^n}\bar \Delta\bigr)$-MMP with scaling of $\rho^*\pi^*H$. Since both $\diff_{\bar D^n}\bar \Delta$ and
$\rho^*\pi^*H$ are $\tau$-invariant, the $\tau$-action descends to
give (regular) involutions 
$$
\tau^i:  \bigl((\bar D^{i})^n, \diff_{\bar D^n}\bar \Delta^{i}\bigr)
\to
 \bigl((\bar D^{i})^n, \diff_{\bar D^n}\bar \Delta^{i}\bigr).
$$
Let $Z^i\subset \bar X^i$ denote the intersection of $D^i$ with the exceptional locus of
$$
(\phi^i\circ\cdots\circ \phi^1)^{-1}:\bar X^i\map \bar X.
$$
 By our assumption, $Z^i$ does not contain any of the log centers of $\bigl(\bar X^i, \bar D^i+\bar\Delta^{i}\bigr) $.
Thus  $\tau^i$ defines a finite equivalence relation on $\bar X^i$ by \cite[9.55]{kk-singbook}.  Therefore the  geometric quotient  $\pi^i: \bar X^i\to X^i$  of $\bar X^i$ by the equivalence relation generated by $\tau^i$ exists by \cite[9.21]{kk-singbook}.   Next \cite[5.38]{kk-singbook} shows that 
 $(X^i, \Delta^i)$ is slc. By Lemma \ref{lem2}  the resulting rational map
$$
\phi^i: (X^{i-1}, \Delta^{i-1})\map (X^i, \Delta^i)
$$
is an MMP step with scaling of $H$.  \qed


\medskip

Note that if  $X$ is  a normal crossing variety \cite[1.7]{kk-singbook} then the log centers of $(X,0)$ 
are exactly the log canonical centers of $(X,0)$, which are also the strata of $X$,
so the important distinction between  log  centers and log canonical centers is not visible in this case.

The normalization  $\pi: (\bar X, \bar D)\to X$ is a normal crossing pair.
It is conjectured that $(\bar X, \bar D)$ has a minimal model. This is currently known if $K_{\bar X}+\bar D$ has non-negative Kodaira dimension  (on every irreducible component) and the dimension is $\leq 5$  \cite{birkar-2007}.

If  a minimal model  $\phi: X\map X^{\rm min}$ exists, then its normalization $(\bar X^{\rm min}, \bar D^{\rm min})$  is a dlt pair whose canonical class is nef and big. The abundance conjecture predicts that 
its canonical class is semi-ample, but  this is known only if   the dimension is $\leq 4$  \cite{MR3567594}.  
However, if abundance holds for $(\bar X^{\rm min}, \bar D^{\rm min})$  then 
 \cite{hacon-xu-brep} implies that the canonical class of $X^{\rm min}$   is also semi-ample.
In particular, the  canonical ring of $X$ is finitely generated.

Thus Theorem \ref{at2}  implies the following. Conjecturally, the dimension restrictions should not be necessary.

\begin{cor} Let $X$ be a pure dimensional, projective, normal crossing variety. Assume that $K_X$ has non-negative Kodaira dimension  on every irreducible component of $X$ and its stable base locus does not contain any stratum of $X$. 
\begin{enumerate} 
\item If  $\dim X\leq 5$  then $X$ has a minimal model  $\phi: X\map X^{\rm min}$,   $\phi$ is a local isomorphism at all log canonical centers and
$X^{\rm min}$  is semi-dlt  \cite[5.19]{kk-singbook}. 
\item If  $\dim X\leq 4$  then the canonical ring of $X$ is finitely generated. \qed
\end{enumerate}
\end{cor}

Before we start the proof of Theorem \ref{at1}, we need to define what a step of an MMP is.

\begin{defn}[MMP steps] \label{say1}
An {\it MMP step} is a diagram of $S$-schemes
$$
\begin{array}{rcl} (X, \Theta) &\stackrel{\phi}{\map} &(X',\Theta')\\ 
f  &\searrow \quad\swarrow & f'\\  &Z&
\end{array}
\eqno{(\ref{say1}.1)}
$$
with the following properties.
\begin{enumerate}\setcounter{enumi}{1}
\item   $(X,\Theta)$ and $(X',\Theta')$ are  pure dimensional lc pairs,
\item $\phi$ is birational, 
\item  $f, f'$ are projective and generically finite,
\item  $-(K_X+\Theta)$ is $f$-ample and 
$K_{X'}+\Theta'$ is $f'$-ample,
\item $f'$ has no exceptional divisors and
\item  $\Theta'=\phi_*\Theta$.
\end{enumerate}
Note that (3) and (6) together imply that $\phi$ is a rational contraction, that is, $\phi^{-1}$ has no exceptional divisors.

For slc pairs, one needs to pay extra attention to the non-normal locus, and there are various possible definitions. However, if $\phi$ is a local isomorphism at all codimension 1 singular points, then the above definition works without changes. This is the only case that we use in the sequel.

We frequently call $\phi: (X, \Theta)\map (X',\Theta')$ an MMP step if
it sits in a diagram as in (\ref{say1}.1) for suitable $Z$. 
Note that $Z$ is not uniquely determined by $\phi: (X, \Theta)\map (X',\Theta')$; if $Z\to Z_1$ is finite then we can replace $Z$ by $Z_1$. 
The usual choice is to take the unique $Z$ such that $f_*\o_X=\o_Z$.
However, the latter condition is not preserved when passing to the normalization of $X$ or to a divisor in $X$. 
Thus allowing different choices of $Z$ is  convenient for us.

If $H$ is a $\q$-Cartier divisor on $X$ then (\ref{say1}.1) is an MMP step {\it with scaling of $H$} if, in addition, 
\begin{enumerate}\setcounter{enumi}{7}
\item  $H$ is $f$-ample, $-H':=-\phi_* H$ is $f'$-ample,
\item $K_X+\Theta+cH$ is numerically $f$-trivial
for some $c\in \q$, (this implies that  $K_{X'}+\Theta'+cH'$ is numerically $f'$-trivial) and
\item $K_X+\Theta+cH$ has positive degree on every proper, irreducible  curve $C\subset X$ that is not contracted by $f$ (and lies over a closed point of $Z$).
\end{enumerate}

In practice we start with a pair  $(X, \Theta+c'H)$ such that
$K_X+\Theta+c'H$ is  ample over $S$. We then decrease the value of $c'$ until we reach $c\leq c'$ such that  $K_X+\Theta+cH$ is nef but not  ample.
If a multiple of  $K_X+\Theta+cH$ is semiample, it gives us
$f:X\to Z$; see \cite{bchm} for details. 
\end{defn}

The following comparison result is clear from the definition.

\begin{lem} \label{lem2} Let   $(X,\Theta)$ and $(X',\Theta')$ be  pure dimensional slc pairs with normalizations 
$\pi: (\bar X, \bar D+\bar \Theta)\to (X, \Theta)$ and
$\pi': (\bar X', \bar D'+\bar \Theta')\to (X', \Theta')$. Then 
 (\ref{say1}.1) is an MMP step iff 
$$
\begin{array}{rcl} (\bar X, \bar D+\bar \Theta) &\stackrel{\bar \phi}{\map} &(\bar X',\bar D'+\bar \Theta')\\ 
\bar f  &\searrow \quad\swarrow & \bar f'\\  &Z&
\end{array}
\eqno{(\ref{lem2}.1)}
$$
is an  MMP step, where $\bar f=f\circ \pi$ and  $\bar f'=f'\circ \pi'$.

Furthermore, if $H$ is a $\q$-Cartier divisor on $X$ then (\ref{say1}.1) is an MMP step  with scaling of $H$ iff (\ref{lem2}.1) is an   MMP step with scaling of $\pi^*H$. \qed
\end{lem}

Next we consider a generalization of MMP steps.

\begin{defn} \label{say3}
A diagram as in (\ref{say1}.1) is called a 
 {\it sub-MMP step} if
\begin{enumerate}
\item  the assumptions  (\ref{say1}.2--5) hold,
\item $f'$ is allowed to have exceptional divisors and
\item $\coeff_{G'}\Theta'\leq \coeff_{G'}\Theta$ for every divisor $G'\subset X'$ that is not $f'$-exceptional. 
(By Lemma \ref{mmp.restricts} this inequality then holds for all divisors over $X$.)
\end{enumerate}
The following example is good to keep in mind.
Let $X$ be a smooth surface and $C\subset X$ a smooth, rational curve with self-intersection $\leq -3$.  Let  $X\to X'$ denote the contraction of $C$. 

Then  $(X,C)\map  (X, 0)$ and  $(X',0)\map  (X, 0)$  are both   sub-MMP step.
Thus $\phi$ can be an isomorphism on the underlying varieties yet a non-trivial sub-MMP step. 
\end{defn}

The main reason for this definition is Lemma \ref{lem4}, but
first we prove that 
the  usual discrepancy inequalities  
(cf.\ \cite[3.38]{km-book} or \cite[1.19 and 1.22]{kk-singbook}) also hold for
sub-MMP steps.

\begin{lem}\label{mmp.restricts} Consider  a sub-MMP step of lc pairs
$$
\begin{array}{rcl} (X, \Theta) &\stackrel{\phi}{\map} &(X',\Theta')\\ 
f  &\searrow \quad\swarrow & f'\\  &Z&
\end{array}
$$
where  $f,f'$ are birational. 
 Then  $a(E, X', \Theta')\geq  a(E, X, \Theta)$
 for every divisor $E$ over $X$. Furthermore,  for every $E$,  the following are equivalent.
\begin{enumerate}
\item $a(E, X', \Theta')>  a(E, X, \Theta)$.
\item $\phi$ is not a local isomorphism at the generic point of $\cent_XE$.
\item $\phi^{-1}$ is not a local isomorphism at the generic point of $\cent_{X'}E$.
\item Either $f$ or $f'$ has positive dimensional fiber  over the generic point of $\cent_ZE$.
\end{enumerate}
\end{lem}

Proof.  
Let $Y$ be the normalization of the  main component of the fiber 
product $X\times_ZX'$ with projections
$X\stackrel{g}{\leftarrow} Y\stackrel{g'}{\to} X'$.
Write
$$
K_Y\simq g^*(K_X+\Theta)-F \qtq{and}K_Y\simq g'^*(K_{X'}+\Theta')-F'
\eqno{(\ref{mmp.restricts}.5)}
$$
where $g_*F=\Theta$ and $g'_*F'=\Theta'$.
Thus
$$
F'-F\simq g'^*(K_{X'}+\Theta') - g^*(K_X+\Theta)
\qtq{is  $(f'\circ g')$-nef.}
\eqno{(\ref{mmp.restricts}.6)}
$$
Note that $(f'\circ g')_*(F-F')=f_*\Theta-f'_*\Theta'$ is 
effective by assumption (\ref{say3}.3). Therefore $F-F'$ is effective by
\cite[3.39]{km-book}, proving the required  inequality.

It is clear that (1) $\Rightarrow$ (2),  (2) $\Leftrightarrow$ (3) and  (2) $\Rightarrow$ (4). Thus assume (4). 

By \cite[3.39]{km-book} the  support of $F-F'$ contains $\ex(f'\circ g')$.
Arguing similarly we get that it also contains $\ex(f\circ g)$. Thus $a(E, X', \Theta')>  a(E, X, \Theta)$
if either  $f$ or $f'$ has positive dimensional fiber  over the generic point of $\cent_ZE$. \qed

\begin{cor}\label{mmp.restricts.cor} A sub-MMP step 
 $\phi: (X, \Theta)\map (X',\Theta')$ is an MMP step iff
  $a(G', X', \Theta')=  a(G', X, \Theta)$
 for every divisor $G'\subset X'$.
\end{cor}

Proof. If $\phi$ is an  MMP step then $\Theta'=\phi_*\Theta$, hence 
  $a(G', X', \Theta')=  a(G', X, \Theta)$  for every divisor $G'\subset X'$.

Conversely, if $G'\subset X'$ is an $f'$-exceptional divisor then
$a(G', X', \Theta')>  a(G', X, \Theta)$ by Lemma \ref{mmp.restricts}.2. Thus there are no $f'$-exceptional divisors and so $\Theta'=\phi_*\Theta$. \qed

\medskip

\begin{lem} \label{lem4} Let $\phi: (X, \Theta)\map (X',\Theta')$ be an MMP step  sitting in a diagram  (\ref{say1}.1). Assume that
 $(X, \Theta)$ is lc,   $\Theta=D+\Delta$ where $D$ is reduced with
normalization $\rho:\bar D\to D$ and 
 none of the irreducible components of $D$ is contracted by $\phi$. 
 Then the  diagram
$$
\begin{array}{rcl} \bigl(\bar D, \diff_{\bar D}\Delta\bigr) &\stackrel{\phi_D}{\map} & \bigl(\bar D', \diff_{\bar D'}\Delta'\bigr)\\ f_D
&\searrow \quad\swarrow & f'_D\\  &Z&
\end{array}
\eqno{(\ref{lem4}.1)}
$$
is a sub-MMP step.
\end{lem}

Proof.  Assumptions  (\ref{say1}.2--4) are  clear and
(\ref{say1}.5) holds since
$$
K_{\bar D}+\diff_{\bar D}\Delta\simq  \rho^*(K_X+D+\Delta).
$$
  It remains to show that (\ref{say3}.3) holds.
More generally, we show that  
$$
a(E, \bar D, \diff_D\Delta)\leq  
 a\bigl(E,\bar D', \diff_{D'}\Delta'\bigr)
\eqno{(\ref{lem4}.2)}
$$ for every divisor $E$ over $\bar D$.

We may assume that $f,f'$ are birational. 
Let $Y$ be the normalization of the  main component of the fiber 
product $X\times_ZX'$ with projections
$X\stackrel{g}{\leftarrow} Y\stackrel{g'}{\to} X'$.
As in (\ref{mmp.restricts}.5) write
$$
g^*(K_X+D+\Delta)\simq g'^*\bigl(K_{X'}+D'+\Delta'\bigr)+F-F',
\eqno{(\ref{lem4}.3)}
$$
where $F-F'$ is effective by \cite[3.38]{km-book}  or by Lemma \ref{mmp.restricts}.

Let $D_Y$ denote the normalization of the birational transform of
$D$ on $Y$. Restricting (\ref{lem4}.3) to $D_Y$ we get
$$
(g|_{D_Y})^*(K_{\bar D}+\diff_{\bar D}\Delta)\simq 
(g'|_{D_Y})^*\bigl(K_{\bar D'}+\diff_{\bar D'}\Delta'\bigr)+F|_{D_Y}
\eqno{(\ref{lem4}.4)}
$$
and $F|_{D_Y}$ is also effective.  \qed

\begin{cor} \label{lem4.cor} Using the notation and assumptions of
Lemma \ref{lem4}, let $p\in \bar D$ be a point. Then $\phi_D: \bigl(\bar D, \diff_{\bar D}\Delta\bigr){\map}  \bigl(\bar D', \diff_{\bar D'}\Delta'\bigr)$ is a 
local isomorphism at $p$ iff $\phi: X\map X'$ is a 
local isomorphism at $\pi(p)$.
\end{cor}

Note that the claims about $X$ and $D$ are different. As in Example \ref{xmp1}.1, 
it can happen that $\phi_D:\bar D\map \bar D'$ is an 
 isomorphism  but $\diff_{\bar D'}\Delta'\neq (\phi_D)_*\diff_{\bar D}\Delta$. 
 
\medskip

Proof. If $\phi$  is a 
local isomorphism at $\pi(p)$ then clearly $\phi_D$ is a 
local isomorphism at $p$. Conversely, if  $\phi_D:\bar D\map \bar D'$ is a 
local isomorphism at $p$ then the maps
$g_D:D_Y\to \bar D$ and $g'_D:D_Y\to \bar D'$ are isomorphic to each other near $p$. 
By (\ref{lem4}.4)
$$
g_D^*\diff_{\bar D}\Delta - {g'_D}^*\diff_{\bar D'}\Delta'=
(g|_{D_Y})^*(F-F').
$$
If $\phi$  is not a 
local isomorphism at $\pi(p)$ then $\supp(F-F')$ contains $p$ by
 by \cite[3.38]{km-book}  or by Lemma \ref{mmp.restricts}, thus    $\diff_{\bar D}\Delta \neq \diff_{\bar D'}\Delta' $
in every neighborhood of $p$. \qed

\begin{prop} \label{prop4}
Using the notation of Lemma \ref{lem4}, assume in addition that 
$D\cap \ex(\phi)$ 
does not contain any log center of  $(X, D+\Delta)$.
Then (\ref{lem4}.1) is an MMP step. 
\end{prop}

Proof.  Assume to the contrary that  (\ref{lem4}.1) is not an an MMP step.
Then, by Corollary \ref{mmp.restricts.cor}, there is a divisor $G'\subset \bar D'$
such that   
$$
a(G', \bar D, \diff_{\bar D}\Delta)< a(G', \bar D', \diff_{\bar D'}\Delta').
\eqno{(\ref{prop4}.1)}
$$ 
Since $a(G', \bar D', \diff_{\bar D'}\Delta')=-\coeff_{G'} \diff_{\bar D'}\Delta'\leq 0$, this implies that
$\cent_{\bar D}G'$ is a log center of
$(\bar D, \diff_{\bar D}\Delta)$. By adjunction \cite[4.8]{kk-singbook},
$\cent_{X}G'$ is also a log center of $(X, D+\Delta)$.

Finally (\ref{prop4}.1) also shows that  $\phi$ is not a  local isomorphism at the generic point of  $\cent_{X}G'$.
\qed

\begin{say}[Proof of Theorem \ref{at1}] \label{pf.at1} By assumption none of the irreducible components of $D$ is contained in $\ex(\Phi^m)$, thus the
maps $\phi^i_{D}$ are birational. They sit in  diagrams
$$
\begin{array}{rcl} \bigl(\bar D^{i-1}, \diff_{\bar D}\Delta^{i-1}\bigr) &\stackrel{\phi^i_D}{\map} & \bigl(\bar D^{i}, \diff_{\bar D}\Delta^{i}\bigr)\\ f^i_D
&\searrow \quad\swarrow & g^i_D\\  &Z^i&
\end{array}
\eqno{(\ref{pf.at1}.1)}
$$
that are sub-MMP steps by Lemma \ref{lem4}. 

If $\phi^m_D$ is not an MMP step then, by Corollary \ref{mmp.restricts.cor}, 
there is a divisor  $G^m\subset \bar D^m$
such that 
$$
a(G^m, \bar D^{m-1},\diff_{\bar D}\Delta^{m-1})< 
a(G^m, \bar D^{m},\diff_{\bar D}\Delta^{m-1})\leq 0.
$$
Combining with the inequalities
$$
a(G^m, \bar D^{i-1},\diff_{\bar D}\Delta^{i-1})\leq
a(G^m, \bar D^{i},\diff_{\bar D}\Delta^{i})
$$
of Lemma \ref{mmp.restricts}.1, we get that
$$
a(G^m, \bar D,\diff_{\bar D}\Delta)< 
a(G^m, X^{m},\diff_{\bar D}\Delta^{m-1})\leq 0.
$$
Thus $\cent_{\bar D}G^m$ is log center of $(\bar D,\diff_{\bar D}\Delta)$.
By adjunction  \cite[4.8]{kk-singbook}, 
its image in $X$ is a  log center of $(X,D+\Delta)$ that is contained in $\ex(\Phi^m)$. \qed
\end{say}






\medskip

Note that  Proposition \ref{prop4} almost implies Theorem \ref{at1}, except that it is not quite clear how to compare $\ex(\Phi^m)\subset X$ with the $\ex(\phi^i)\subset X^{i-1}$ that are needed to directly apply  Proposition \ref{prop4}. The following variant of the concept of  exceptional set
gives a clearer picture and a slighty different way of deriving 
Theorem \ref{at1}. 

\begin{defn}[Divisorial exceptional set]\label{divex.def}
Let $\phi:X\map X'$ be a birational map of schemes that are proper over $S$. 
The {\it divisorial exceptional set} of $\phi$, denoted by $\dex(\phi)$, is 
the set of all divisors $E$ over $X$ such that
$\phi$ is not a local isomorphism at the generic point of $\cent_XE$.

Thus the usual exceptional set $\ex(\phi)\subset X$ is the union of the centers of the divisors in $\dex(\phi)$. The advantage of
 divisorial exceptional sets is that we can compare them for
different birational models.

\end{defn}

\begin{lem} \label{divex.lem}
Let  $\phi_i:(X^{i-1}, \Delta^{i-1})\map (X^i, \Delta^i)$ be a sequence of MMP steps. Then
\begin{enumerate}
\item $\dex(\phi^m\circ\cdots\circ \phi^1)=\{E: a(E, X^0, \Delta^0) < a(E, X^m, \Delta^m)\}$ and
\item $\dex(\phi^m\circ\cdots\circ \phi^1)=\dex(\phi^1)\cup\cdots\cup\dex(\phi^m)$.
\end{enumerate}
\end{lem}

Proof.  The containments
$$
\begin{array}{rcl}
\dex(\phi^m\circ\cdots\circ \phi^1) & \supset &\{E: a(E, X^0, \Delta^0) < a(E, X^m, \Delta^m)\}\\
\dex(\phi^m\circ\cdots\circ \phi^1)& \subset &\dex(\phi^1)\cup\cdots\cup\dex(\phi^m)
\end{array}
\eqno{(\ref{divex.lem}.3)}
$$
are clear.
 For a single MMP step $\phi:(X, \Delta)\map (X', \Delta')$, 
\cite[3.38]{km-book} shows that 
$$
\dex(\phi)=\{E: a(E, X, \Delta) < a(E, X', \Delta')\}.
\eqno{(\ref{divex.lem}.4)}
$$
Combining with the inequalities 
$a(E,X^{i-1}, \Delta^{i-1})\leq a(E,X^i, \Delta^i)$ we obtain that
$a(E,X^{0}, \Delta^{0})\leq a(E,X^m, \Delta^m)$ and 
$$
a(E,X^{0}, \Delta^{0})< a(E,X^m, \Delta^m) \Leftrightarrow
E\in \dex(\phi^1)\cup\cdots\cup\dex(\phi^m).
$$
This shows that
$$
\{E: a(E,X^{0}, \Delta^{0})< a(E,X^m, \Delta^m)\}= 
\dex(\phi^1)\cup\cdots\cup\dex(\phi^m),
$$
which completes the proof. \qed

\def\cprime{$'$} \def\cprime{$'$} \def\cprime{$'$} \def\cprime{$'$}
  \def\cprime{$'$} \def\cprime{$'$} \def\cprime{$'$} \def\dbar{\leavevmode\hbox
  to 0pt{\hskip.2ex \accent"16\hss}d} \def\cprime{$'$} \def\cprime{$'$}
  \def\polhk#1{\setbox0=\hbox{#1}{\ooalign{\hidewidth
  \lower1.5ex\hbox{`}\hidewidth\crcr\unhbox0}}} \def\cprime{$'$}
  \def\cprime{$'$} \def\cprime{$'$} \def\cprime{$'$}
  \def\polhk#1{\setbox0=\hbox{#1}{\ooalign{\hidewidth
  \lower1.5ex\hbox{`}\hidewidth\crcr\unhbox0}}} \def\cdprime{$''$}
  \def\cprime{$'$} \def\cprime{$'$} \def\cprime{$'$} \def\cprime{$'$}
\providecommand{\bysame}{\leavevmode\hbox to3em{\hrulefill}\thinspace}
\providecommand{\MR}{\relax\ifhmode\unskip\space\fi MR }
\providecommand{\MRhref}[2]{%
  \href{http://www.ams.org/mathscinet-getitem?mr=#1}{#2}
}
\providecommand{\href}[2]{#2}

\bigskip

\noindent Institute of Mathematics ``Simion Stoilow'' of the Romanian Academy,

\noindent  P.O.\ BOX 1-764, RO-014700 Bucharest, Romania.

{\begin{verbatim} florin.ambro@imar.ro\end{verbatim}}
\medskip

\noindent  Princeton University, Princeton NJ 08544-1000

{\begin{verbatim} kollar@math.princeton.edu\end{verbatim}}

\end{document}